\documentstyle[12pt]{article}
\begin{document}

\title{Fixed point free homeomorphisms of the complex plane}

\author{Nikolaos E. Sofronidis\footnote{$A \Sigma MA:$ 130/2543/94}}

\date{\footnotesize Department of Economics, University of Ioannina, Ioannina 45110, Greece.
(nsofron@otenet.gr, nsofron@cc.uoi.gr)}

\maketitle

\begin{abstract}
Our purpose in this article is to prove that the group $H({\bf C})$ of homeomorphisms of the complex plane ${\bf C}$ is a metric group equipped with the metric induced by uniform convergence of homeomorphisms and their inverses on compacts and the set $$\left\{ h \in H({\bf C}) : ( \forall z \in {\bf C} )( h(z) \neq z ) \right\} $$ of fixed point free homeomorphisms of the complex plane is a conjugacy invariant dense $G_{\delta }$ subset of $H({\bf C})$.
\end{abstract}

\section*{\footnotesize{{\bf Mathematics Subject Classification:} 03E15, 37E30, 57S05.}}

{\bf 1. Definition.} Following page 808 of [4], if $$d_{u}(f,g) = \sum\limits_{n=1}^{ \infty } 2^{-n} \frac{ \sup\limits_{ \vert z \vert \leq n } \vert f(z) - g(z) \vert }{ 1 + \sup\limits_{ \vert z \vert \leq n } \vert f(z) - g(z) \vert }$$ whenever $f$, $g$ are continuous functions ${\bf C} \rightarrow {\bf C}$, then following 8) on page 60 of [2], it is not difficult to verify that if $d(f,g) = d_{u}(f,g) + d_{u} \left( f^{-1} , g^{-1} \right) $, whenever $f$, $g$ are homeomorphisms of ${\bf C}$, then the group $H({\bf C})$ of homeomorphisms of ${\bf C}$ equipped with the metric $d$ becomes a metric space and convergence in $d$ is equivalent to uniform convergence of homeomorphisms and their inverses on compacts, so that for any compact subset $K$ of ${\bf C}$ and for any open subset $O$ of ${\bf C}$, the set $\left\{ f \in H({\bf C}) : f[K] \subseteq O \right\} $ is open in $H({\bf C})$ equipped with $d$, as it follows from Problem 8 b. on page 193 of [3].
\\ \rm \\
{\bf 2. Definition.} If $X$ is any compact Polish space, then we denote by $K(X)$ the compact Polish space of compact subsets of $X$ equipped with the Hausdorff metric and for any sequence $\left( K_{n} \right) _{n\in {\bf N}}$ of compact subsets of $X$, we denote by $\liminf\limits_{n \rightarrow \infty }K_{n}$ the topological lower limit of $\left( K_{n} \right) _{n\in {\bf N}}$ and we denote by $\limsup\limits_{n \rightarrow \infty }K_{n}$ the topological upper limit of $\left( K_{n} \right) _{n\in {\bf N}}$, while if these two limits coincide, then we call them the topological limit $\lim\limits_{n \rightarrow \infty }K_{n}$ of $\left( K_{n} \right) _{n\in {\bf N}}$, which is equivalent to convergence in $K(X)$ with respect to the Hausdorff metric on it. See, for example, Section 4.F on pages 24-28 of [2].
\\ \rm \\
{\bf 3. Lemma.} If $h_{n} \rightarrow h$ in $H({\bf C})$ as $n \rightarrow \infty $, then for any $K \in K({\bf C})$, we have that $\lim\limits_{n \rightarrow \infty }h_{n}[K]=h[K]$ in $K({\bf C})$.
\\ \rm \\
{\bf Proof.} If $p \in K$, then evidently $h(p) = \lim\limits_{n \rightarrow \infty }h_{n}(p) \in \liminf\limits_{n \rightarrow \infty }h_{n}[K]$ and consequently $h[K] \subseteq \liminf\limits_{n \rightarrow \infty }h_{n}[K]$. So let $q \in \limsup\limits_{n \rightarrow \infty }h_{n}[K]$ be arbitrary but fixed. For any $k \in {\bf N}$, there exists an integer $n_{k} > n_{k-1}$, where $n_{-1} = -1$, such that $h_{n_{k}}[K] \cap D \left( q ; 2^{-k} \right) \neq \emptyset $ and hence there exists $p_{k} \in K$ for which $h_{n_{k}} \left( p_{k} \right) \in D \left( q ; 2^{-k} \right) $. The compactness of $K$ implies that there exists a subsequence $\left( p_{k_{j}} \right) _{j \in {\bf N}}$ of $\left( p_{k} \right) _{k \in {\bf N}}$ which converges to some point $p \in K$. Hence, since $\left( h_{n} \right) _{n \in {\bf N}}$ converges to $h$ uniformly on each compact subset of ${\bf C}$, by virtue of Exercise 40 e. on page 162 of [3], it follows that $\left( h_{n_{k_{j}}} \right) _{j \in {\bf N}}$ converges continu- ously to $h$ and consequently $q = \lim\limits_{j \rightarrow \infty }h_{n_{k_{j}}} \left( p_{k_{j}} \right) = h(p) \in h[K]$. Therefore, we have that $\limsup\limits_{n \rightarrow \infty }h_{n}[K] \subseteq h[K]$ and consequently $\lim\limits_{n \rightarrow \infty }h_{n}[K]=h[K]$, thus the claim follows. \hfill $\bigtriangleup $
\\ \rm \\
{\bf 4. Lemma.} If $g_{n} \rightarrow g$ and $h_{n} \rightarrow h$ in $H({\bf C})$ as $n \rightarrow \infty $, then $\left( g_{n} \circ h_{n} \right) \rightarrow (g \circ h)$ in $H({\bf C})$ as $n \rightarrow \infty $.
\\ \rm \\
{\bf Proof.} Let $K$ be any non-empty compact subset of ${\bf C}$. Then ${\overline{B}}(h[K];1)$ is compact in ${\bf C}$. Indeed, the compactness of $h[K]$, as it follows from Proposition 24 of [3], implies that there exist $q_{1}$, ..., $q_{k}$ in $h[K]$, where $k$ is a positive integer, for which $h[K] \subseteq \bigcup\limits_{j=1}^{k} D \left( q_{j} ; 1 \right) $ and let $p \in {\overline{B}}(h[K];1)$ be arbitrary but fixed. Then the compactness of $h[K]$ implies that there exists $q \in h[K]$ such that $\vert p-q \vert = \inf\limits_{z \in h[K]} \vert p-z \vert \leq 1$ and hence there exists $j \in \{ 1, ..., k \}$ for which $q \in D \left( q_{j} ; 1 \right) $, therefore $\left\vert p - q_{j} \right\vert \leq \vert p - q \vert + \left\vert q - q_{j} \right\vert \leq 1 + 1 = 2$, i.e., $p \in {\overline{D}} \left( q_{j} ; 2 \right) $. We have thus proved that ${\overline{B}}(h[K];1) \subseteq \bigcup\limits_{j=1}^{k} {\overline{D}} \left( q_{j} ; 2 \right) $ and the claim follows from the fact that the ${\overline{D}} \left( q_{j} ; 2 \right) $'s are compact. See, for example, Proposition 22 on page 156 of [3]. By virtue of Lemma 3 and Section 4.F on pages 24-28 of [2], it follows that there exists $N \in {\bf N}$ such that for any integer $n \geq N$, we have that $h_{n}[K] \subseteq B(h[K];1)$. Therefore, given any integer $n \geq N$, we have that
\begin{enumerate}
\item[ ]
$\sup\limits_{p \in K} \left\vert \left( g_{n} \circ h_{n} \right) (p) - (g \circ h)(p) \right\vert $
\item[ ]
$\leq \sup\limits_{p \in K} \left\vert g_{n} \left( h_{n} (p) \right) - g \left( h_{n} (p) \right) \right\vert  + \sup\limits_{p \in K} \left\vert g \left( h_{n} (p) \right) - g(h(p)) \right\vert $
\item[ ]
$\leq \sup\limits_{q \in {\overline{B}}(h[K];1)} \left\vert g_{n}(q) - g(q) \right\vert + \sup\limits_{p \in K} d_{g} \left( h_{n}(p), h(p) \right) $,
\end{enumerate}
where $d_{g}(p,q) = \vert g(p)-g(q) \vert $, whenever $p$, $q$ range over ${\bf C}$, also constitutes a complete compatible metric on ${\bf C}$, so $\left( g_{n} \circ h_{n} \right) _{n \in {\bf N}}$ converges uniformly to $(g \circ h)$ on compacts. An analogous argument shows that $\left( h_{n}^{-1} \circ g_{n}^{-1} \right) _{n \in {\bf N}}$ converges uniformly to $\left( h^{-1} \circ g^{-1} \right) $ on compacts, so the claim follows. \hfill $\bigtriangleup $
\\ \rm \\
{\bf 5. Corollary.} $H({\bf C})$ equipped with the metric $d$ is a metric group.
\\ \rm \\
{\bf 6. Definition.} If $u(z) = \frac{z}{1 - \vert z \vert }$, whenever $z \in D(0;1)$, then it is not difficult to verify that $u$ is a homeomorphism $D(0;1) \rightarrow {\bf C}$ whose inverse is defined by the relation $u^{-1}(w) = \frac{w}{ 1 + \vert w \vert }$, whenever $w \in {\bf C}$, so if $H(D(0;1))$ is the group of homeomorphisms of $D(0;1)$, then $H(D(0;1)) \ni \psi \mapsto \left( u \circ \psi \circ u^{-1} \right) \in H({\bf C})$ is a group isomorphism via which the metric group $H({\bf C})$ makes $H(D(0;1))$ also a metric group in which convergence is also uniform convergence on compacts.
\\ \rm \\
{\bf 7. Definition.} If $F$ is any non-empty subset of ${\bf C}$, then we call $F$ a closed $2$-cell in ${\bf C}$, if there exist a closed disk ${\overline{D}}( \alpha ; \rho )$ contained in $D(0;1)$ and a homeomorphism $k$ of $D(0;1)$ onto an open subset $U$ of ${\bf C}$ for which $F = k \left[ {\overline{D}}( \alpha ; \rho ) \right] $. In addition, if $h \in H({\bf C})$, then we set $${\bf C} \setminus supp(h) = \bigcup\limits_{O \in {\bf \Sigma }_{1}^{0}({\bf C})} \left\{ z \in O : h(z)=z \right\} $$ where ${\bf \Sigma }_{1}^{0}({\bf C})$ is the set of open subsets of ${\bf C}$. See, for example, 11.B on page 68 of [2].
\\ \rm \\
{\bf 8. Theorem.} If $F$ is any closed $2$-cell in ${\bf C}$, then $\left\{ h \in H({\bf C}) : supp(h) \subseteq F \right\} $ is closed nowhere dense in $H({\bf C})$.
\\ \rm \\
{\bf Proof.} Let ${\overline{D}}( \alpha ; \rho )$ be the closed disk contained in $D(0;1)$ and let $k$ be the homeomorphism of $D(0;1)$ onto the open subset $U$ of ${\bf C}$ for which $F = k \left[ {\overline{D}}( \alpha ; \rho ) \right] $. We will first prove that $\left\{ h \in H({\bf C}) : supp(h) \subseteq F \right\} $ is closed in $H({\bf C})$. Let $\eta > 0$ be such that ${\overline{D}}( \alpha ; \rho + 2 \eta ) \subseteq D(0;1)$. Then, evidently
\begin{enumerate}
\item[ ]
$\left\{ h \in H({\bf C}) : supp(h) \subseteq F \right\} $
\item[ ]
$= \bigcap\limits_{ 0 < \delta \leq \eta } \left\{ h \in H({\bf C}) : supp(h) \subseteq k \left[ D( \alpha ; \rho + 2 \delta ) \right] \right\} $.
\end{enumerate}
Hence, if $h_{n} \rightarrow h$ in $H({\bf C})$ as $n \rightarrow \infty $ and for any $n \in {\bf N}$, we have that $supp \left( h_{n} \right) \subseteq F$, then for any $\delta \in (0, \eta ]$ and for any $n \in {\bf N}$, we have that $h_{n} = id$ in ${\bf C} \setminus k \left[ D( \alpha ; \rho + 2 \delta ) \right] $, which is obviously closed in ${\bf C}$, and consequently $h = id$ in ${\bf C} \setminus k \left[ D( \alpha ; \rho + 2 \delta ) \right] $, whenever $\delta \in (0, \eta ]$, which is easily seen to imply that $h=id$ in $F$ and the claim follows. Thus, given any $h \in H({\bf C})$ for which $supp(h) \subseteq F$, what is left to show is that $h$ is not an interior point of $\left\{ g \in H({\bf C}) : supp(g) \subseteq F \right\} $. Given any $\delta \in [0, \eta ]$, we set
\[
{\psi }_{ \delta } \left( \alpha + r e^{i \theta } \right) = \left\{
\begin{array}{lllll}
\alpha + \frac{ \rho + \delta }{ \rho } r e^{i \theta } & \mbox{if $0 \leq r \leq \rho $}
\\
 & \mbox{ and $0 \leq \theta < 2 \pi $}
\\ \\
\alpha + \left( \frac{1}{2} (r - \rho ) + \rho + \delta \right) e^{i \theta } & \mbox{if $\rho \leq r \leq \rho + 2 \delta $}
\\
 & \mbox{ and $0 \leq \theta < 2 \pi $}
\\ \\
\alpha + r e^{i \theta } & \mbox{otherwise}
\end{array}
\right.
\]
whenever $r \geq 0$ and $0 \leq \theta < 2 \pi $ are such that $\left( \alpha + r e^{i \theta } \right) \in D(0;1)$. It is not difficult to see that ${\psi }_{\delta } : D(0;1) \rightarrow D(0;1)$ constitutes a homeomorphism which expands radially ${\overline{D}}( \alpha ; \rho )$ to ${\overline{D}}( \alpha ; \rho + \delta )$ and shrinks radially $\left( {\overline{D}}( \alpha ; \rho + 2 \delta ) \setminus D( \alpha ; \rho ) \right) $ to $\left( {\overline{D}}( \alpha ; \rho + 2 \delta ) \setminus D( \alpha ; \rho + \delta ) \right) $, while ${\psi }_{ \delta } = id$ in $\left( D(0;1) \setminus {\overline{D}}( \alpha ; \rho + 2 \delta ) \right) $. Moreover, it is not difficult to verify that for any $\epsilon $, $\delta $ in $[0, \eta ]$ and for any $z \in D(0;1)$, we have that $\left\vert {\psi }_{\epsilon }(z) - {\psi }_{\delta }(z) \right\vert \leq \vert \epsilon - \delta \vert $, so $\sup\limits_{z \in D(0;1)}\left\vert {\psi }_{\epsilon }(z) - {\psi }_{\delta }(z) \right\vert \leq \vert \epsilon - \delta \vert $. It is not difficult to see that ${\psi }_{\delta }^{-1} : D(0;1) \rightarrow D(0;1)$ constitutes a homeomorphism which shrinks radially ${\overline{D}}( \alpha ; \rho + \delta )$ to ${\overline{D}}( \alpha ; \rho )$ and expands radially $\left( {\overline{D}}( \alpha ; \rho + 2 \delta ) \setminus D( \alpha ; \rho + \delta ) \right) $ to $\left( {\overline{D}}( \alpha ; \rho + 2 \delta ) \setminus D( \alpha ; \rho ) \right) $, while ${\psi }_{ \delta }^{-1} = id$ in $\left( D(0;1) \setminus {\overline{D}}( \alpha ; \rho + 2 \delta ) \right) $, which implies that $\left\vert {\psi }_{\epsilon }^{-1}(z) - {\psi }_{\delta }^{-1}(z) \right\vert \leq \vert \epsilon - \delta \vert $, so $\sup\limits_{z \in D(0;1)}\left\vert {\psi }_{\epsilon }^{-1}(z) - {\psi }_{\delta }^{-1}(z) \right\vert \leq \vert \epsilon - \delta \vert $. Therefore, the mapping $[0, \eta ] \ni \delta \mapsto {\psi }_{\delta } \in H(D(0;1))$ is continuous. Given any $\delta \in [0, \eta ]$, let $h_{\delta } : {\bf C} \rightarrow {\bf C}$ be defined by the relation
\[
h_{\delta }(p) = \left\{
\begin{array}{ll}
\left( k \circ {\psi }_{\delta } \circ k^{-1} \right) (p) & \mbox{if $p \in k \left[ {\overline{D}}( \alpha ; \rho + 2 \delta ) \right] $}
\\ \\
p & \mbox{otherwise}
\end{array}
\right.
\]
whenever $p \in {\bf C}$. Then, it is not difficult to see that $h_{\delta } : {\bf C} \rightarrow {\bf C}$ constitutes a homeomorphism which expands $k \left[ {\overline{D}}( \alpha ; \rho ) \right] $ to $k \left[ {\overline{D}}( \alpha ; \rho + \delta ) \right] $ and shrinks $\left( k \left[ {\overline{D}}( \alpha ; \rho + 2 \delta ) \right] \setminus k \left[ D( \alpha ; \rho ) \right] \right) $ to $\left( k \left[ {\overline{D}}( \alpha ; \rho + 2 \delta ) \right] \setminus k \left[ D( \alpha ; \rho + \delta ) \right] \right) $, while $h_{\delta } = id$ in $\left( {\bf C} \setminus k \left[ {\overline{D}}( \alpha ; \rho + 2 \delta ) \right] \right) $. Moreover, the continuity of the mapping $[0, \eta ] \ni \delta \mapsto {\psi }_{\delta } \in H(D(0;1))$ is easily seen to imply that $[0, \eta ] \ni \delta \mapsto h_{\delta } \in H({\bf C})$ is continuous and hence so is $[0, \eta ] \ni \delta \mapsto \left( h_{\delta } \circ h \right) \in H({\bf C})$, as it follows from Lemma 4. Thus, the claim follows from the fact that $h_{0} = id$ and for any $\delta \in (0, \eta ]$, we have that $supp \left( h_{\delta } \circ h \right) $ is not a subset of $k \left[ {\overline{D}}( \alpha ; \rho ) \right] $, since $supp(h) \subseteq k \left[ {\overline{D}}( \alpha ; \rho ) \right] $ and $h_{\delta }$ expands $k \left[ {\overline{D}}( \alpha ; \rho ) \right] $ to $k \left[ {\overline{D}}( \alpha ; \rho + \delta ) \right] $. \hfill $\bigtriangleup $
\\ \rm \\
{\bf 9. Lemma.} For any $c \in {\bf C}$, we have that $\left\{ h \in H({\bf C}) : h(c)=c \right\} $ is closed nowhere dense in $H({\bf C})$.
\\ \rm \\
{\bf Proof.} It is not difficult to see that $\left\{ h \in H({\bf C}) : h(c)=c \right\} $ is closed in $H({\bf C})$ and let $g \in \left\{ h \in H({\bf C}) : h(c)=c \right\} $. If we set ${\tau }_{a}(z) = z+a$, whenever $z$, $a$ range in ${\bf C}$, then since for any $a$, $b$ in ${\bf C}$, we have that $\sup\limits_{z \in {\bf C}} \left\vert {\tau }_{a}(z) - {\tau }_{b}(z) \right\vert = \vert a-b \vert $, it follows immediately that ${\bf C} \ni a \mapsto {\tau }_{a} \in H({\bf C})$ constitutes a continu- ous and injective homomorphism between the corresponding metric groups and consequently the mapping ${\bf C} \ni a \mapsto \left( {\tau }_{a} \circ g \right) \in H({\bf C})$ is continuous, as it follows from Lemma 4, which implies that $\left( {\tau }_{a} \circ g \right) \rightarrow g $ in $H({\bf C})$ as $a \rightarrow 0$ in ${\bf C}$. Thus, the claim follows from the fact that for any $a \in \left( {\bf C} \setminus \{ 0 \} \right) $, we have that $\left( {\tau }_{a} \circ g \right) (c) = c+a \neq c$. \hfill $\bigtriangleup $
\\ \rm \\
{\bf 10. Theorem.} $\left\{ h \in H({\bf C}) : ( \forall z \in {\bf C} )( h(z) \neq z ) \right\} $ constitutes a conjugacy invariant dense $G_{\delta }$ subset of $H({\bf C})$.
\\ \rm \\
{\bf Proof.} If $C$ is any countable dense subset of ${\bf C}$, then by virtue of Lemma 9, it is enough to show that $$\left\{ h \in H({\bf C}) : ( \forall z \in {\bf C} )( h(z) \neq z ) \right\} = \bigcap\limits_{c \in C} \left\{ h \in H({\bf C}) : h(c) \neq c \right\} $$ or (equivalently) that $$\bigcap\limits_{c \in C} \left\{ h \in H({\bf C}) : h(c) \neq c \right\} \subseteq \left\{ h \in H({\bf C}) : ( \forall z \in {\bf C} )( h(z) \neq z ) \right\} $$ Indeed, if $g \in \bigcap\limits_{c \in C} \left\{ h \in H({\bf C}) : h(c) \neq c \right\} $, then it is not difficult to prove that for any $c \in C$, there exists ${\epsilon }_{c} > 0$ such that $g \left[ D \left( c ; {\epsilon }_{c} \right) \right] \cap D \left( c ; {\epsilon }_{c} \right) = \emptyset $ and since $\bigcup\limits_{c \in C} D \left( c ; {\epsilon }_{c} \right) = {\bf C}$, it follows immediately that $$g \in \left\{ h \in H({\bf C}) : ( \forall z \in {\bf C} )( h(z) \neq z ) \right\} $$ \hfill $\bigtriangleup $
\\ \rm \\
{\bf 11. Remark.} The reader is referred to [1] and its references for other results regarding fixed point free homeomorphisms of the complex plane.

\end{document}